\documentclass[12pt,twoside]{article}
\usepackage{amssymb}
\usepackage{amsmath}

\date{}

\begin{document}
\centerline {\Large{ Characterizing of Inner Product Spaces  by the Mapping $n_{x,y}$ }}

\centerline{}

\centerline{}
\centerline{ Hossein Dehghan
%\footnote{   }
}
\centerline{}
\centerline{Department of Mathematics, Institute for
Advanced Studies in Basic
 Sciences
 }
\centerline{ (IASBS), Gava Zang,  Zanjan 45137-66731, Iran}
\centerline{ Email:  hossein.dehgan@gmail.com}
\centerline{}

\newtheorem{theorem}{Theorem}[section]
\newtheorem{lemma}[theorem]{Lemma}
\newtheorem{proposition}[theorem]{Proposition}
\newtheorem{corollary}[theorem]{Corollary}
\newtheorem{definition}[theorem]{Definition}
\newtheorem{example}[theorem]{Example}
\newtheorem{xca}[theorem]{Exercise}
\newtheorem{remark}[theorem]{Remark}
\newtheorem{algorithm}{Algorithm}
\numberwithin{equation}{section}

\begin{abstract}
For the vectors $x$ and $y$ in a normed linear spaces $X$, the mapping $n_{x,y}: \mathbb{R}\to \mathbb{R}$ is defined by  $n_{x,y}(t)=\|x+ty\|$. In this note, comparing the mappings $n_{x,y}$ and $n_{y,x}$ we obtain a simple and useful characterization of inner product spaces.
 
 %To do this we first show that it we first study  the notion of skew-angular distance in normed linear spaces. Then, we study the relation between skew-angular distance and the mapping $n_{x,y}(t)=\|x+ty\|$ in normed linear spaces. Finally, we obtain a simple and useful characterization of inner product spaces.

\end{abstract}


\section{ Introduction}
\par
Let $(X,\|.\|)$ be a real normed linear space and $x,y\in X$.
In 1998, Dragomir and Koliha \cite{DK,DK2} introduced the following mapping and studied some its properties.
\begin{eqnarray}
\nonumber n_{x,y}: \mathbb{R} \rightarrow\mathbb{R}, \ \ \ n_{x,y}(t)=\|x+ty\|,
\end{eqnarray}
where $\mathbb{R}$ denotes the set of real numbers. They proved that $n_{x,y}$ is convex, continuous and has one sided derivatives at each point of $\mathbb{R}$ \cite[Proposition 2.1]{DK}.
In 2006, Maligranda \cite[Theorem 1]{Mal1} (see also \cite{Mal2}) introduced the following strengthening of the triangle inequality and its reverse: For any nonzero vectors $x$ and $y$ in a real normed linear space $X = (X, \|.\|)$ it is true that
\begin{eqnarray}\label{Ma re tri up}
\|x+y\|\leq\|x\|+\|y\|-\left(2-\left\|\frac{x}{\|x\|}+\frac{y}{\|y\|}\right\|\right)
\min\{\|x\|,\|y\|\}
\end{eqnarray}
and
\begin{eqnarray}\label{Ma re tri low}
\|x+y\|\geq\|x\|+\|y\|-\left(2-\left\|\frac{x}{\|x\|}+\frac{y}{\|y\|}\right\|\right)
\max\{\|x\|,\|y\|\}.
\end{eqnarray}
%Moreover, if either $\|x\|=\|y\|$ or $y=cx$ with $c>0$, then equality holds in both (\ref{Ma re tri up}) and (\ref{Ma re tri low}).
Also, the author used  (\ref{Ma re tri up}) and (\ref{Ma re tri low}) for the following estimation of the \emph{angular
distance}  $\alpha[x,y]=\|\frac{x}{\|x\|}-\frac{y}{\|y\|}\|$ between two nonzero elements $x$ and $y$ in $X$ which was defined by Clarkson  in \cite{Clarc}.
\begin{eqnarray}\label{angu}
\frac{\|x-y\|-|\ \|x\|-\|y\|\ |}{\min\{\|x\|,\|y\|\}}\leq\alpha[x,y]\leq\frac{\|x-y\|+|\ \|x\|-\|y\|\ |}{\max\{\|x\|,\|y\|\}}.
\end{eqnarray}
The right hand of estimate (\ref{angu}) is a refinement of the Massera-Schaffer inequality proved in 1958
(see \cite[ Lemma 5.1]{MS}): for nonzero
vectors $x$ and $y$ in $X$ we have that $\alpha[x,y]\leq\frac{2\|x-y\|}{\max\{\|x\|,\|y\|\}}$ , which is stronger than
the Dunkl-Williams inequality $\alpha[x,y]\leq\frac{4\|x-y\|}{\|x\|+\|y\|}$ proved in \cite{DW}. In the same paper, Dunkl and Williams proved that the constant 4 can be replaced by 2 if and only if $X$ is an inner product space. 
On the other hand, the author \cite{Deh} introduced skew-angular distance as follows.
\begin{definition}
 For two nonzero elements $x$ and $y$ in a real normed linear space $X=(X,\|.\|)$  the distance
\begin{eqnarray}
\beta[x,y]=\left\|\frac{x}{\|y\|}-\frac{y}{\|x\|}\right\|
\end{eqnarray}
is called skew-angular distance between  $x$ and $y$.
\end{definition}
In the same paper, a characterizations of inner product spaces obtained by comparing angular and skew-angular distance.
\section{Main results}
In this section, we first compare the functions $n_{x,y}$ and $n_{y,x}$ with each other and we show this is equivalent to compare the angular distance $\alpha[x,y]$ with the skew-angular distance $\beta[x,y]$. Using these results we present a simple characterization of inner product spaces.
\begin{proposition}
Let $(X,\|.\|)$ be a normed linear space. Then the following statements are mutually equivalent:
\begin{enumerate}
  \item[(i)] $n_{x,y}(t)\leq n_{y,x}(t)$ for all $x,y\in X$ and all $t\in [0,1]$.
  \item[(ii)]  $n_{x,y}(t)\leq n_{y,x}(t)$ for all $x,y\in X$ and all $t\in [-1,0]$.
  \item[(iii)] $n_{y,x}(t)\leq n_{x,y}(t)$ for all $x,y\in X$ and all $t\in [1,\infty)$.
    \item[(iv)] $n_{y,x}(t)\leq n_{x,y}(t)$ for all $x,y\in X$ and all $t\in (-\infty,-1]$.
\end{enumerate}
\end{proposition}
\textbf{ Proof.} Since $n_{x,y}(t)=n_{x,-y}(-t)$ and $n_{y,x}(t)=n_{-y,x}(-t)$, then $(i)\Leftrightarrow (ii)$ and $(iii)\Leftrightarrow (iv)$. If $t\not= 0$, then
$$\|y+tx\|\leq \|x+ty\|\ \ \ \Leftrightarrow\ \ \ \left\|x+\frac{1}{t}y\right\|\leq\left\|y+\frac{1}{t}x\right\|.$$
Moreover, continuity of $n_{x,y}$ implies that
 $$n_{x,y}(0)=\lim_{n\to\infty}n_{x,y}(1/n)=\lim_{n\to\infty}n_{x,y}(-1/n).$$
Thus $(i)\Leftrightarrow (iii)$ and $(ii)\Leftrightarrow (iv)$. This completes proof. $\Box$\\
\par
By the above proposition, it is sufficient to compare the functions $n_{x,y}$ and $n_{y,x}$ on the interval $[0,1]$.
\begin{lemma}\label{ineq in inner}
For any vectors $x$ and $y$ in an inner product space  $(X,\langle.,.\rangle)$, it is true that
\begin{eqnarray}
n_{x^*,y^*}(t)\leq n_{y^*,x^*}(t)
\end{eqnarray}
for all $t\in [0,1]$, where $x^*, y^*$ is rearrangement of $x,y$ with $\|x^*\|\leq\|y^*\|$.
\end{lemma}
\textbf{ Proof.} It follows from polar identity that
\begin{eqnarray}
\nonumber n^2_{x^*,y^*}(t)-n^2_{y^*,x^*}(t)=\left(\|y^*\|^2-\|x^*\|^2\right)t^2-\left(\|y^*\|^2-\|x^*\|^2\right)\leq0
\end{eqnarray}
for all $t\in[0,1]$. $\Box$\\
The following example shows that Lemma \ref{ineq in inner} is not true in general normed spaces.
\begin{example}  Let $X=\mathbb{R}^2$ with the norm of $x=(a,b)$ given by $\|x\|=|a|+|b|$. Taking $x=(0,2)$ and $y=(2,-1)$, we see that
 $\|x\|\leq\|y\|$ but $n_{x,y}(1/2)> n_{y,x}(1/2)$.
\end{example}

 The next theorem due to Lorch will be useful in the sequel.
\begin{theorem}\label{Lorch th}
(See \cite{Lorch}.) Let $(X,\|.\|)$ be a real normed linear space. Then the following statements are mutually equivalent:
\begin{enumerate}
  \item[(i)] For each $x,y\in X$ if $\|x\|=\|y\|$, then $\|x+y\|\leq \|\gamma x+\gamma^{-1}y\|$ (for all $\gamma\neq 0$).
  \item[(ii)] For each $x,y\in X$ if $\|x+y\|\leq \|\gamma x+\gamma^{-1}y\|$ (for all $\gamma\neq 0$), then $\|x\|=\|y\|$.
  \item[(iii)] $(X,\|.\|)$  is an inner product space.
  \end{enumerate}
\end{theorem}
\par
Now, we state and prove the main theorem which seems new and simple.
\begin{theorem}
Let $(X,\|.\|)$ be a normed linear space.
%and $x^*, y^*$ denotes rearrangement of $x,y$ in $X$ with $\|x^*\|\leq\|y^*\|$.
 Then the following statements are mutually equivalent:
\begin{enumerate}
  \item[(i)] $n_{x,y}(t)\leq n_{y,x}(t)$ for all $x,y\in X$ with $\|x\|\leq\|y\|$ and all $t\in [0,1]$.
  \item[(ii)]  $\alpha[x,y]\leq\beta[x,y]$ for all nonzero vectors $x,y\in X$.
\item[(iii)]  $(X,\|.\|)$  is an inner product space.
\end{enumerate}
\end{theorem}
\textbf{ Proof.} Let $x,y\in X$ and $x,y\not=0$. Without loss of generality we may assume that $\|x\|\leq\|y\|$. Using $(i)$ for $t=\|x\|/\|y\|$ and  $x,-y$ instead of $x,y$ we obtain that
$$\left\|x-\frac{\|x\|}{\|y\|}y\right\|\leq \left\|y-\frac{\|x\|}{\|y\|}x\right\|.$$
Since $x\not=0$, then
$$\left\|\frac{x}{\|x\|}-\frac{y}{\|y\|}\right\|\leq \left\|\frac{x}{\|y\|}-\frac{y}{\|x\|}\right\|.$$
This proves $(i)\Rightarrow (ii)$.\\
 $(ii)\Rightarrow (iii)$ follows from \cite[Theorem3.2]{Deh} but for convenience of the reader we inset it here. 
 %Next, we show that $(ii)$ implies $(i)$ of Theorem \ref{Lorch th}.
  Let $x,y\in X$, $\|x\|=\|y\|$ and $\gamma\neq 0$. From Theorem \ref{Lorch th} it is enough to prove that  $\|x+y\|\leq \|\gamma x+\gamma^{-1}y\|$. If $x=0$ or $y=0$, then the proof is clear. Let $x\neq 0$, $y\neq 0$ and $\gamma>0$. Applying inequality $\alpha[x,y]\leq\beta[x,y]$  to  $\gamma^{\frac{1}{2}} x$ and $-\gamma^{-\frac{1}{2}}y$ instead of $x$ and $y$, respectively, we obtain
\begin{eqnarray}\label{}\nonumber
\left\|\frac{\gamma^{\frac{1}{2}}x}{\gamma^{\frac{1}{2}}\|x\|}+\frac{\gamma^{-\frac{1}{2}}y}{\gamma^{-\frac{1}{2}}\|y\|}\right\|
\leq\left\|\frac{\gamma^{\frac{1}{2}}x}{\gamma^{-\frac{1}{2}}\|y\|}+\frac{\gamma^{-\frac{1}{2}}y}{\gamma^{\frac{1}{2}}\|x\|}\right\|.
\end{eqnarray}
Thus
\begin{eqnarray}\label{}\nonumber
\left\|\frac{x}{\|x\|}+\frac{y}{\|y\|}\right\|
\leq\left\|\frac{\gamma x}{\|y\|}+\frac{\gamma^{-1}y}{\|x\|}\right\|.
\end{eqnarray}
Since $\|x\|=\|y\|\neq0$, then
$$\|x+y\|\leq \|\gamma x+\gamma^{-1}y\|.$$
Now, let $\gamma$ be negative. Put $\mu=-\gamma>0$. From the positive case we get
$$\|x+y\|\leq  \|\mu x+\mu^{-1}y\|=\|\gamma x+\gamma^{-1}y\|.$$
$(iii)\Rightarrow (i)$ follows from Lemma \ref{ineq in inner}.
This completes the proof. $\Box$\\

\textbf{Acknowledgment}\\
The author thanks  Professor J. Rooin for his valuable suggestions which improved the original manuscript.

\end{document}